\documentclass[11pt,twoside,reqno]{article}
\usepackage{fullpage}
\usepackage{amssymb,amsfonts,amsmath,amsthm,bm}

\usepackage{ stmaryrd}
\theoremstyle{definition}
	\newtheorem{thm}{Theorem}

	\newtheorem{conj}[thm]{Conjecture}
	
\usepackage{enumerate}
\usepackage{enumitem}
\usepackage{hyperref}
\usepackage{textcomp}
\usepackage{color}

\usepackage[mathscr]{eucal}
\usepackage{bbm}

\usepackage{tikz}

\usepackage{hyperref}
\usepackage[dvipsnames]{xcolor}
\newcommand\myshade{85}
\colorlet{mylinkcolor}{red}
\colorlet{mycitecolor}{blue}
\colorlet{myurlcolor}{Aquamarine}
\hypersetup{
	linkcolor  = mylinkcolor,
	citecolor  = mycitecolor,
	urlcolor   = myurlcolor!\myshade!black,
	colorlinks = true,
}

\newcommand{\bi}{\mathbf{i}}
\newcommand{\bj}{\mathbf{j}}
\newcommand{\bk}{\mathbf{k}}

\newcommand{\rel}{\text{rel}}

\newcommand{\relh}{{\widehat{\rel}}}

\newcommand{\ds}{\displaystyle}

\usepackage[dvipsnames]{xcolor}

\newcommand{\GG}{{\textcolor{Green}{\bf G}}}
\newcommand{\RR}{{\textcolor{Red}{\bf R}}}
\newcommand{\BB}{{\textcolor{Blue}{\bf B}}}

\newcommand{\qtri}[5]{\begin{bmatrix} #1 \\ #2, \, #3, \, #4 \end{bmatrix}_{#5}}

\numberwithin{equation}{section}
\allowdisplaybreaks

\newcommand{\ignore}[1]{}

\usepackage{array}

\begin{document}
\thispagestyle{empty}
\title{On a pair of three-colored (mod 10) partition identities}

\author{Matthew C. Russell
\footnote{Department of Statistics, Texas A\&M University \url{matthewcrussell@tamu.edu}.}}

\maketitle

\begin{abstract}
We prove a pair of (mod 10) partition identities. The sum sides involve three-colored partitions into distinct parts, while the product sides are the generating functions for distinct partitions times the Rogers-Ramanujan products. Our proofs make heavy use of Maple to verify that functional equations are satisfied.
\end{abstract}

\section{Introduction}
A partition of $n$ is a finite non-increasing sequence $\left(\lambda_1,\dots,\lambda_k\right)$ of positive integers (called {\it parts}) that sums to $n$: $\lambda_1+\cdots+\lambda_k=n.$ We can consider colored partitions by choosing some set of colors, assigning a color to each part, and considering parts that are the same integer but different colors to be distinct.

Let $\Gamma$ be the set of partitions into three colors ($R$, $G$, and $B$), subject to the following conditions:
\begin{itemize}
\item The sizes of all parts must be distinct. (For all $j$, there is at most one total copy of $j_\RR$, $j_\GG$, and $j_\BB$.)
\item The following subpartitions are forbidden for all $j$:
\begin{itemize}
\item[$\star$] $j_\RR + \left(j+1\right)_\RR$
\item[$\star$] $j_\RR + \left(j+2\right)_\RR$
\item[$\star$] $j_\RR + \left(j+1\right)_\BB$
\item[$\star$] $j_\GG + \left(j+1\right)_\RR$
\item[$\star$] $j_\GG + \left(j+2\right)_\RR$
\item[$\star$] $j_\GG + \left(j+1\right)_\BB$
\item[$\star$] $j_\BB + \left(j+1\right)_\RR$
\item[$\star$] $j_\BB + \left(j+1\right)_\GG$
\end{itemize}
\end{itemize}
Here is a graphical representation of the difference conditions of $\Gamma$. If a line connects two parts, then that pair of parts is forbidden from appearing together in a colored partition in $\Gamma$:
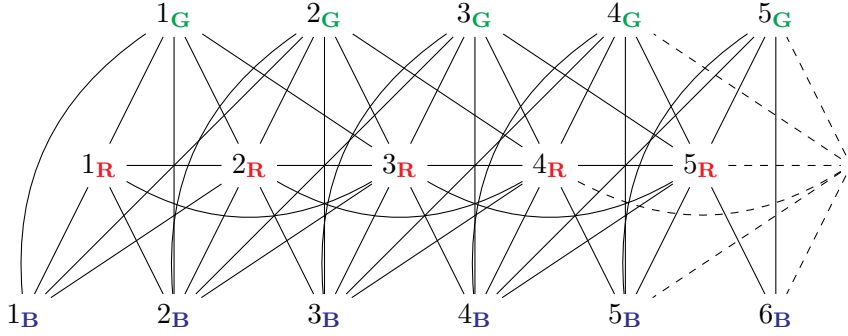
\begin{center}
\begin{tikzpicture}
        \node (c1r) at (0,0) {$1_\RR$};
        \node (c2r) at (2,0) {$2_\RR$};
        \node (c3r) at (4,0) {$3_\RR$};
        \node (c4r) at (6,0) {$4_\RR$};
        \node (c5r) at (8,0) {$5_\RR$};
        \node (c6r) at (10,0) {};
        \node (c7r) at (12,0) {};
        \node (c1g) at (1,2) {$1_\GG$};
        \node (c2g) at (3,2) {$2_\GG$};
        \node (c3g) at (5,2) {$3_\GG$};
        \node (c4g) at (7,2) {$4_\GG$};
        \node (c5g) at (9,2) {$5_\GG$};
        \node (c6g) at (11,2) {};
        \node (c1b) at (-1,-2) {$1_\BB$};
        \node (c2b) at (1,-2) {$2_\BB$};
        \node (c3b) at (3,-2) {$3_\BB$};
        \node (c4b) at (5,-2) {$4_\BB$};
        \node (c5b) at (7,-2) {$5_\BB$};
        \node (c6b) at (9,-2) {$6_\BB$};
        \draw[-] (c1r) -- (c1g) ;
        \draw[-] (c2r) -- (c2g) ;
        \draw[-] (c3r) -- (c3g) ;
        \draw[-] (c4r) -- (c4g) ;
        \draw[-] (c5r) -- (c5g) ;
        \draw[-] (c1r) -- (c1b) ;
        \draw[-] (c2r) -- (c2b) ;
        \draw[-] (c3r) -- (c3b) ;
        \draw[-] (c4r) -- (c4b) ;
        \draw[-] (c5r) -- (c5b) ;
        \draw[dashed] (c6r) -- (c6b) ;
        \draw[-] (c1g) to[bend right] (c1b) ;
        \draw[-] (c2g) to[bend right] (c2b) ;
        \draw[-] (c3g) to[bend right] (c3b) ;
        \draw[-] (c4g) to[bend right] (c4b) ;
        \draw[-] (c5g) to[bend right] (c5b) ;
        \draw[-] (c2r) -- (c1b) ;
        \draw[-] (c3r) -- (c2b) ;
        \draw[-] (c4r) -- (c3b) ;
        \draw[-] (c5r) -- (c4b) ;
        \draw[dashed] (c6r) -- (c5b) ;
        \draw[-] (c2r) -- (c1g) ;
        \draw[-] (c3r) -- (c2g) ;
        \draw[-] (c4r) -- (c3g) ;
        \draw[-] (c5r) -- (c4g) ;
        \draw[dashed] (c6r) -- (c5g) ;
        \draw[-] (c1r) -- (c2b) ;
        \draw[-] (c2r) -- (c3b) ;
        \draw[-] (c3r) -- (c4b) ;
        \draw[-] (c4r) -- (c5b) ;
        \draw[-] (c5r) -- (c6b) ;
        \draw[-] (c3r) -- (c1g) ;
        \draw[-] (c4r) -- (c2g) ;
        \draw[-] (c5r) -- (c3g) ;
        \draw[dashed] (c6r) -- (c4g) ;
        \draw[-] (c1b) --  (c2g);
        \draw[-] (c2b) --  (c3g);
        \draw[-] (c3b) --  (c4g);
        \draw[-] (c4b) --  (c5g);
        \draw[-] (c1g) --  (c2b);
        \draw[-] (c2g) --  (c3b);
        \draw[-] (c3g) --  (c4b);
        \draw[-] (c4g) --  (c5b);
        \draw[-] (c5g) --  (c6b);
        \draw[-] (c1r) -- (c2r);
        \draw[-] (c2r) -- (c3r);
        \draw[-] (c3r) -- (c4r);
        \draw[-] (c4r) -- (c5r);
        \draw[dashed] (c5r) -- (c6r);
        \draw[-] (c1r) to[bend right] (c3r);
        \draw[-] (c2r) to[bend right] (c4r);
        \draw[-] (c3r) to[bend right] (c5r);
        \draw[dashed] (c4r) to[bend right] (c6r);
\end{tikzpicture}        
\end{center}
Alternatively, we can express the difference conditions by the following matrix (as was done, for example, by K. Alladi, G. E. Andrews, and B. Gordon~\cite{AAG}):
\begin{center}
\begin{tabular}{c|ccc}
     & $\RR$ & $\GG$  & $\BB$ \\ \hline
$\RR$ & 3 & 1 & 2\\
$\GG$ & 3 & 1 & 2 \\
$\BB$ & 2 & 2 & 1\\
\end{tabular}
\end{center}
The way to read this is by rows. For example, if there is a part $j_\RR$, for example, then the next largest possible red part would be $(j+3)_\RR$, the next largest possible green part would be $(j+1)_\GG$, and the next largest possible blue part would be $(j+2)_\BB$.

Throughout this paper, $q$ will be a formal variable. We use the standard $q$-Pochhammer notation $(a;q)_n = \prod_{j=1}^n \left(1-aq^{j-1}\right)$ for $j \in \{0,1,2,\dots\}\cup\infty$, along with $\left(a_1,a_2,\dots,a_k;q\right)_n=\left(a_1;q\right)_n\left(a_2;q\right)_n\cdots\left(a_k;q\right)_n$.
The goal of this paper is to prove the following pair of identities involving three-colored partitions in $\Gamma$:
\begin{thm} \label{thm:main1} Let $A(n)$ be the number of three-colored partitions in $\Gamma$ that sum to a non-negative integer $n$ with no occurrences of $1_\RR$. Then, 
\begin{equation} \ds\sum_{n \ge 0} A(n) q^n = \frac{1}{\left(q;q^2\right)_\infty\left(q^1,q^4;q^5\right)_\infty}. \label{eq:main1}
\end{equation}
\end{thm}

\begin{thm} \label{thm:main2} Let $A^\ast(n)$ be the number of three-colored partitions in $\Gamma$ that sum to a non-negative integer $n$ with no occurrences of $1_\RR$, $2_\RR$, and $1_\BB$. Then, 
\begin{equation} \ds\sum_{n \ge 0} A^\ast(n) q^n = \frac{1}{\left(q;q^2\right)_\infty\left(q^2,q^3;q^5\right)_\infty}. \label{eq:main2}
\end{equation} 
\end{thm}

To illustrate the above identities, the $A(6)=18$ colored partitions in $\Gamma$ with no occurrences of $1_\RR$ (satisfying the conditions of Theorem~\ref{thm:main1}) are 
\begin{align*}
&6_\RR     &&6_\GG    && 6_\BB    && 5_\RR+1_\GG&& 5_\RR+1_\BB   && 5_\GG+1_\GG \\
&5_\GG+1_\BB &&5_\BB+1_\GG&& 5_\BB+1_\BB&& 4_\RR+2_\BB&& 4_\GG+2_\RR   && 4_\GG+2_\GG \\
&4_\GG+2_\BB &&4_\BB+2_\RR&&4_\BB+2_\GG && 4_\BB+2_\BB&&3_\GG+2_\GG+1_\GG&&3_\BB+2_\BB+1_\BB,
\end{align*}
and the $A^\ast(6)=12$ colored partitions in $\Gamma$ with no occurrences of $1_\RR$, $2_\RR$, and $1_\BB$ (satisfying the conditions of Theorem~\ref{thm:main2}) are 
\begin{align*}
&6_\RR&&6_\GG&&6_\BB&&5_\RR+1_\GG&& 5_\GG+1_\GG&& 5_\BB+1_\GG\\
&4_\RR+2_\BB&&4_\GG+2_\GG&&4_\GG+2_\BB&&4_\BB+2_\GG&&4_\BB+2_\BB&&3_\GG+2_\GG+1_\GG.
\end{align*}
It is easily checked that the coefficients in the Taylor series expansions of the right sides of~\eqref{eq:main1} and~\eqref{eq:main2} are 18 and 12, respectively.

These identities were discovered through a computerized search for colored partition identities similar to the one previously carried out by S.~Kanade and the author~\cite{KR1}. Other, similar identities were found in a related manner by Kanade, D. Nandi, and the author~\cite{KNR}. These products are also the principally specialized characters of level 3 standard modules of $A_1^{(1)}$. Other papers with identities exploiting this connection include works by A. Meurman and M. Primc~\cite{MeuPri-99}, S. Chern, Z. Li, D. Stanton, T. Xue, and A.~J. Yee~\cite{CheEtAl}, J. Dousse, L. Hardiman, and I. Konan~\cite{DouHarKon}, and Kanade and the author~\cite{KR5}; J. Stembridge also found an identity with the product side of Theorem~\ref{thm:main2}~\cite{Stem}. Another family of colored partition identities in which all parts must be distinct was found by Dousse~\cite{DouSch}. Many papers have been written generalizing Capparelli's identities~\cite{Capp} to colored partitions~\cite{AAG,BM,DouPrimc}. Dousse~\cite{D_Sil} found a colored partition generalization and refinement of a partition identity of I. Siladi\'c~\cite{Sil}; Dousse and J. Lovejoy~\cite{DL} found a colored partition generalization and refinement of a partition identity of Primc~\cite{Primc}. Many other papers have dealt with identities for colored partitions; we briefly mention three more recent ones~\cite{AndColor,Zad,MCR_Mac}.

\subsection{Acknowledgments.} The author thanks Shashank Kanade for many continued conversations and for writing the Maple implementation of the Murray-Miller algorithm that was adapted for this paper, and also thanks Premila Samuel Russell for graciously providing feedback on an initial draft of this paper that improved the exposition.

\section{Preliminaries}
We define the $q$-trinomial coefficient $\qtri{i+j+k}ijkq$ as
\begin{align*}
\qtri{i+j+k}ijkq &= \frac{\left(q;q\right)_{i+j+k}}{\left(q;q\right)_i\left(q;q\right)_j\left(q;q\right)_k}.
\end{align*}
Note that this is a different $q$-trinomial coefficient than the one that has frequently appeared in the literature~\cite{AB,AndCapp}. Our $q$-trinomial coefficients satisfy the following analogues of the $q$-Pascal triangle relations:
\begin{align} 
\label{eq:qtriPasc1}
\qtri{i+j+k}ijkq&=&q^{j+k}&\qtri{i+j+k-1}{i-1}jkq+&q^k&\qtri{i+j+k-1}i{j-1}kq+&&\qtri{i+j+k-1}ij{k-1}q \\ 
\label{eq:qtriPasc2}
\qtri{i+j+k}ijkq&=&q^{j+k}&\qtri{i+j+k-1}{i-1}jkq+&&\qtri{i+j+k-1}i{j-1}kq+&q^j&\qtri{i+j+k-1}ij{k-1}q \\
\qtri{i+j+k}ijkq&=&q^k&\qtri{i+j+k-1}{i-1}jkq+&q^{i+k}&\qtri{i+j+k-1}i{j-1}kq+&&\qtri{i+j+k-1}ij{k-1}q \\
\qtri{i+j+k}ijkq&=&&\qtri{i+j+k-1}{i-1}jkq+&q^{i+k}&\qtri{i+j+k-1}i{j-1}kq+&q^i&\qtri{i+j+k-1}ij{k-1}q\\
\qtri{i+j+k}ijkq&=&q^j&\qtri{i+j+k-1}{i-1}jkq+&&\qtri{i+j+k-1}i{j-1}kq+&q^{i+j}&\qtri{i+j+k-1}ij{k-1}q\\ \label{eq:qtriPasc6}
\qtri{i+j+k}ijkq&=&&\qtri{i+j+k-1}{i-1}jkq+&q^i&\qtri{i+j+k-1}i{j-1}kq+&q^{i+j}&\qtri{i+j+k-1}ij{k-1}q
\end{align}

Define, for nonnegative integers $a,b,c$, 
\begin{align*}% \label{eq:SxTx}
S_{a,b,c}(x)&=\sum_{i,j,k\ge 0}
\frac{q^{\frac 32i^2+\frac12j^2+\frac12k^2+ij+ik+jk-\frac12i+\frac12j+\frac12k+ai+bj+ck}}{\left(q;q\right)_{i+j+k}}\qtri{i+j+k}ijk{q^2} x^{2i+j+k} \\
T_{a,b,c}(x)&=\sum_{i,j,k\ge 0}
\frac{q^{\frac 32i^2+\frac12j^2+\frac12k^2+ij+ik+jk-\frac12i+\frac12j+\frac12k+ai+bj+ck}}{\left(q;q\right)_{i+j+k}}\qtri{i+j+k}ijk{q^2} x^{i+j+k}
\end{align*}
Again, $x$ will be another formal variable. We will frequently suppress the argument of $x$ and simply write $S_{a,b,c}$ or $T_{a,b,c}$. Note that $S_{a,b,c}(x)$ and $T_{a,b,c}(x)$ are nearly the same, other than a small change in the exponent of $x$. In particular, observe that 
\begin{equation} \label{eq:29}
S_{a,b,c}(1)=T_{a,b,c}(1).
\end{equation}
Also, note the following shift equations:
\begin{align}
S_{a,b,c}\left(xq^i\right) &= S_{a+2i,b+i,c+i}(x) \label{eq:Sabcshift} \\
T_{a,b,c}\left(xq^i\right) &= T_{a+i,b+i,c+i}(x) \label{eq:Tabcshift}.
\end{align}

\subsection{Atomic relations}
We will use the method of atomic relations, as developed by the author in joint work with Kanade and Baker and Sadowski~\cite{KR4,BKRS}; similar ideas are also present in previous work of Chern~\cite{ChernLinked} (see section 4). Our atomic relations for $S_{a,b,c}$ are:
\begin{align*}
&\rel^1_{a,b,c}: S_{a,b,c}&-&S_{a+2,b,c}&-&x^2q^{a+1}S_{a+3,b+1,c+1}&-&x^2q^{a+2}S_{a+4,b+2,c+2} &&&=0 \\
&\rel^2_{a,b,c}: S_{a,b,c}&-&S_{a,b+2,c}&-&xq^{b+1}S_{a+1,b+3,c+1}&-&xq^{b+2}S_{a+2,b+4,c+2} &&&=0 \\
&\rel^3_{a,b,c}: S_{a,b,c}&-&S_{a,b,c+2}&-&xq^{c+1}S_{a+1,b+1,c+3}&-&xq^{c+2}S_{a+2,b+2,c+4}  &&&=0 \\
&\rel^4_{a,b,c}: S_{a,b,c}&-&S_{a+1,b+1,c+1}&-&x^2q^{a+1}S_{a+3,b+3,c+3}&-&xq^{b+1}S_{a+1,b+1,c+3}&-&xq^{c+1}S_{a+1,b+1,c+1} &=0 \\
&\rel^5_{a,b,c}: S_{a,b,c}&-&S_{a+1,b+1,c+1}&-&x^2q^{a+1}S_{a+3,b+3,c+3}&-&xq^{b+1}S_{a+1,b+1,c+1}&-&xq^{c+1}S_{a+1,b+3,c+1} &=0 \\
&\rel^6_{a,b,c}: S_{a,b,c}&-&S_{a+1,b+1,c+1}&-&x^2q^{a+1}S_{a+3,b+1,c+3}&-&xq^{b+1}S_{a+3,b+1,c+3}&-&xq^{c+1}S_{a+1,b+1,c+1} &=0 \\
&\rel^7_{a,b,c}: S_{a,b,c}&-&S_{a+1,b+1,c+1}&-&x^2q^{a+1}S_{a+3,b+1,c+1}&-&xq^{b+1}S_{a+3,b+1,c+3}&-& xq^{c+1}S_{a+3,b+1,c+1}&=0 \\
&\rel^8_{a,b,c}: S_{a,b,c}&-&S_{a+1,b+1,c+1}&-&x^2q^{a+1}S_{a+3,b+3,c+1}&-&xq^{b+1}S_{a+1,b+1,c+1}&-&xq^{c+1}S_{a+3,b+3,c+1} &=0 \\
&\rel^9_{a,b,c}: S_{a,b,c}&-&S_{a+1,b+1,c+1}&-&x^2q^{a+1}S_{a+3,b+1,c+1}&-&xq^{b+1}S_{a+3,b+1,c+1}&-&xq^{c+1}S_{a+3,b+3,c+1} &=0 
\end{align*}
Here is a proof of $\rel^1_{a,b,c}$. Note the re\"indexing $i-1=\bf{i}$ that occurs between the third and fourth lines, along with the fact that $\ds\frac{\left(q^2;q^2\right)_{\bi+j+k+1}}{\left(q;q\right)_{\bi+j+k+1}} = \frac{\left(q^2;q^2\right)_{\bi+j+k}}{\left(q;q\right)_{\bi+j+k}}\frac{\left(1-q^{2(\bi+j+k+1)}\right)}{\left(1-q^{\bi+j+k+1}\right)}=\frac{\left(1+q^{\bi+j+k+1}\right)\left(q^2;q^2\right)_{\bi+j+k}}{\left(q;q\right)_{\bi+j+k}}$, which is used to go from the antepenultimate line to the penultimate line:
\begin{align*}
 &S_{a,b,c}-S_{a+2,b,c}= S_{a,b,c}\left(1-q^{2i}\right) \\&=\sum_{i,j,k\ge 0}
\frac{q^{\frac 32i^2+\frac12j^2+\frac12k^2+ij+ik+jk-\frac12i+\frac12j+\frac12k+ai+bj+ck}}{\left(q;q\right)_{i+j+k}}\qtri{i+j+k}ijk{q^2} x^{2i+j+k} \left(1-q^{2i}\right) 
\\&=\sum_{i,j,k\ge 0}
\frac{q^{\frac 32i^2+\frac12j^2+\frac12k^2+ij+ik+jk-\frac12i+\frac12j+\frac12k+ai+bj+ck}\left(q^2;q^2\right)_{i+j+k}}{\left(q;q\right)_{i+j+k}\left(q^2;q^2\right)_{i-1}\left(q^2;q^2\right)_{j}\left(q^2;q^2\right)_{k}} x^{2i+j+k}
\\&=\sum_{\bi,j,k\ge 0}
\frac{q^{\frac 32(\bi+1)^2+\frac12j^2+\frac12k^2+(\bi+1)j+(\bi+1)k+jk-\frac12(\bi+1)+\frac12j+\frac12k+a(\bi+1)+bj+ck}\left(q^2;q^2\right)_{\bi+j+k+1}}{\left(q;q\right)_{\bi+j+k+1}\left(q^2;q^2\right)_{\bi}\left(q^2;q^2\right)_{j}\left(q^2;q^2\right)_{k}} x^{2\bi+j+k+2} \\
&=x^2q^{a+1}\sum_{\bi,j,k\ge 0}
\frac{q^{\frac 32\bi^2+\frac12j^2+\frac12k^2+\bi j+\bi k+jk-\frac12\bi+\frac12j+\frac12k+(a+3)\bi+(b+1)j+(c+1)k}\left(q^2;q^2\right)_{\bi+j+k+1}}{\left(q;q\right)_{\bi+j+k+1}\left(q^2;q^2\right)_{\bi}\left(q^2;q^2\right)_{j}\left(q^2;q^2\right)_{k}} x^{2\bi+j+k} \\
&=x^2q^{a+1}\!\!\!\!\sum_{\bi,j,k\ge 0}\!\!\!\!
\frac{q^{\frac 32\bi^2+\frac12j^2+\frac12k^2+\bi j+\bi k+jk-\frac12\bi+\frac12j+\frac12k+(a+3)\bi+(b+1)j+(c+1)k}\left(1+q^{\bi+j+k+1}\right)}{\left(q;q\right)_{\bi+j+k}}\!\!\qtri{\bi+j+k}\bi jk{q^2}  \!\!\!\!x^{2\bi+j+k} \\
&=x^2q^{a+1}S_{a+3,b+1,c+1}+x^2q^{a+2}S_{a+4,b+2,c+2}.
\end{align*}
$\rel^2_{a,b,c}$ and $\rel^3_{a,b,c}$ are proved similarly, involving multiplication by $\left(1-q^{2j}\right)$ and $\left(1-q^{2k}\right)$, respectively.

Now, here is a proof of $\rel^4_{a,b,c}$ using~\eqref{eq:qtriPasc1} (going from the third line to the fourth, fifth, and sixth lines). Three triple sums are produced, and in each one, we shift $i-1={\bf i}$, $j-1={\bf j}$, and $k-1={\bf k}$, respectively:
\begin{align*}
&S_{a,b,c}-S_{a+1,b+1,c+1}= S_{a,b,c}\left(1-q^{i+j+k}\right) \\&=\sum_{i,j,k\ge 0}
\frac{q^{\frac 32i^2+\frac12j^2+\frac12k^2+ij+ik+jk-\frac12i+\frac12j+\frac12k+ai+bj+ck}}{\left(q;q\right)_{i+j+k}}\qtri{i+j+k}ijk{q^2} x^{2i+j+k} \left(1-q^{i+j+k}\right) \\
&=\sum_{i,j,k\ge 0}
\frac{q^{\frac 32i^2+\frac12j^2+\frac12k^2+ij+ik+jk-\frac12i+\frac12j+\frac12k+ai+bj+ck}}{\left(q;q\right)_{i+j+k-1}}\qtri{i+j+k}ijk{q^2} x^{2i+j+k} \\
&=\sum_{i,j,k\ge 0}
\frac{q^{\frac 32i^2+\frac12j^2+\frac12k^2+ij+ik+jk-\frac12i+\frac12j+\frac12k+ai+bj+ck}}{\left(q;q\right)_{i+j+k-1}}q^{2j+2k}\qtri{i+j+k-1}{i-1}jk{q^2} x^{2i+j+k}  \\
&\,\,+\sum_{i,j,k\ge 0}
\frac{q^{\frac 32i^2+\frac12j^2+\frac12k^2+ij+ik+jk-\frac12i+\frac12j+\frac12k+ai+bj+ck}}{\left(q;q\right)_{i+j+k-1}}q^{2k}\qtri{i+j+k-1}i{j-1}k{q^2}x^{2i+j+k} \\
&\,\,+\sum_{i,j,k\ge 0}
\frac{q^{\frac 32i^2+\frac12j^2+\frac12k^2+ij+ik+jk-\frac12i+\frac12j+\frac12k+ai+bj+ck}}{\left(q;q\right)_{i+j+k-1}}\qtri{i+j+k-1}ij{k-1}{q^2} x^{2i+j+k} \\
&=\sum_{\bi,j,k\ge 0}
\frac{q^{\frac 32\bi^2+\frac12j^2+\frac12k^2+\bi j+\bi k+jk-\frac12\bi+\frac12j+\frac12k+a\bi+bj+ck+3\bi+j+k+a+1}}{\left(q;q\right)_{\bi+j+k}}q^{2j+2k}\qtri{\bi+j+k}{\bi}jk{q^2} x^{2\bi+j+k+2}  \\
&\,\,+\sum_{i,\bj,k\ge 0}
\frac{q^{\frac 32i^2+\frac12\bj^2+\frac12k^2+i\bj+ik+\bj k-\frac12i+\frac12\bj+\frac12k+ai+b\bj+ck+i+\bj+k+b+1}}{\left(q;q\right)_{i+\bj+k}}q^{2k}\qtri{i+\bj+k}i{\bj}k{q^2}x^{2i+\bj+k+1} \\
&\,\,+\sum_{i,j,\bk\ge 0}
\frac{q^{\frac 32i^2+\frac12j^2+\frac12\bk^2+ij+i\bk+j\bk-\frac12i+\frac12j+\frac12\bk+ai+bj+c\bk+i+j+\bk+c+1}}{\left(q;q\right)_{i+j+\bk}}\qtri{i+j+\bk}ij{\bk}{q^2} x^{2i+j+\bk+1} \\
&=x^2q^{a+1}\sum_{\bi,j,k\ge 0}
\frac{q^{\frac 32\bi^2+\frac12j^2+\frac12k^2+\bi j+\bi k+jk-\frac12\bi+\frac12j+\frac12k+(a+3)\bi+(b+3)j+(c+3)k}}{\left(q;q\right)_{\bi+j+k}}\qtri{\bi+j+k}{\bi}jk{q^2} x^{2\bi+j+k}  \\
&\,\,+xq^{b+1}\sum_{i,\bj,k\ge 0}
\frac{q^{\frac 32i^2+\frac12\bj^2+\frac12k^2+i\bj+ik+\bj k-\frac12i+\frac12\bj+\frac12k+(a+1)i+(b+1)\bj+(c+3)k}}{\left(q;q\right)_{i+\bj+k}}\qtri{i+\bj+k}i{\bj}k{q^2}x^{2i+\bj+k} \\
&\,\,+xq^{c+1}\sum_{i,j,\bk\ge 0}
\frac{q^{\frac 32i^2+\frac12j^2+\frac12\bk^2+ij+i\bk+j\bk-\frac12i+\frac12j+\frac12\bk+(a+1)i+(b+1)j+(c+1)\bk}}{\left(q;q\right)_{i+j+\bk}}\qtri{i+j+\bk}ij{\bk}{q^2} x^{2i+j+\bk} \\
&=x^2q^{a+1}S_{a+3,b+3,c+3}
+xq^{b+1}S_{a+1,b+1,c+3}+xq^{c+1}S_{a+1,b+1,c+1}.
\end{align*}
The remainder of the atomic relations have similar proofs, using~\eqref{eq:qtriPasc2}--\eqref{eq:qtriPasc6}.

Because of the similarity of $T_{a,b,c}(x)$ to $S_{a,b,c}(x)$, the atomic relations for $T_{a,b,c}$ are nearly the same (other than some factors of $x^2$ being changed to $x$). We omit their proofs.
\begin{align*}
&\relh^1_{a,b,c}: T_{a,b,c}&-&T_{a+2,b,c}&-&xq^{a+1}T_{a+3,b+1,c+1}&-&x^2q^{a+2}T_{a+4,b+2,c+2} &&&=0 \\
&\relh^2_{a,b,c}: T_{a,b,c}&-&T_{a,b+2,c}&-&xq^{b+1}T_{a+1,b+3,c+1}&-&xq^{b+2}T_{a+2,b+4,c+2} &&&=0 \\
&\relh^3_{a,b,c}: T_{a,b,c}&-&T_{a,b,c+2}&-&xq^{c+1}T_{a+1,b+1,c+3}&-&xq^{c+2}T_{a+2,b+2,c+4}  &&&=0 \\
&\relh^4_{a,b,c}: T_{a,b,c}&-&T_{a+1,b+1,c+1}&-&xq^{a+1}T_{a+3,b+3,c+3}&-&xq^{b+1}T_{a+1,b+1,c+3}&-&xq^{c+1}T_{a+1,b+1,c+1} &=0 \\
&\relh^5_{a,b,c}: T_{a,b,c}&-&T_{a+1,b+1,c+1}&-&xq^{a+1}T_{a+3,b+3,c+3}&-&xq^{b+1}T_{a+1,b+1,c+1}&-&xq^{c+1}T_{a+1,b+3,c+1} &=0 \\
&\relh^6_{a,b,c}: T_{a,b,c}&-&T_{a+1,b+1,c+1}&-&xq^{a+1}T_{a+3,b+1,c+3}&-&xq^{b+1}T_{a+3,b+1,c+3}&-&xq^{c+1}T_{a+1,b+1,c+1} &=0 \\
&\relh^7_{a,b,c}: T_{a,b,c}&-&T_{a+1,b+1,c+1}&-&xq^{a+1}T_{a+3,b+1,c+1}&-&xq^{b+1}T_{a+3,b+1,c+3} &-&xq^{c+1}T_{a+3,b+1,c+1}&=0 \\
&\relh^8_{a,b,c}: T_{a,b,c}&-&T_{a+1,b+1,c+1}&-&xq^{a+1}T_{a+3,b+3,c+1}&-&xq^{b+1}T_{a+1,b+1,c+1}&-&xq^{c+1}T_{a+3,b+3,c+1} &=0 \\
&\relh^9_{a,b,c}: T_{a,b,c}&-&T_{a+1,b+1,c+1}&-&xq^{a+1}T_{a+3,b+1,c+1}&-&xq^{b+1}T_{a+3,b+1,c+1}&-&xq^{c+1}T_{a+3,b+3,c+1} &=0
\end{align*}

\subsection{More on $\Gamma$}

For non-negative integers $m$ and $n$: 
\begin{itemize} 
\item Let $A_1(m,n)$ be the number of three-colored partitions in $\Gamma$ of $n$ with exactly $m$ parts that have no occurrences of $1_\RR$.
\item Let $A_2(m,n)$ be the number of three-colored partitions in $\Gamma$ of $n$ with exactly $m$ parts that have no occurrences of $1_\RR$ and $1_\GG$. 
\item Let $A_3(m,n)$ be the number of three-colored partitions in $\Gamma$ of $n$ with exactly $m$ parts that have no occurrences of $1_\RR$, $2_\RR$, and $1_\BB$. 
\end{itemize}
Then, define the generating functions for these three sets as
\begin{align*}
Q_1(x) &= \sum_{m,n \ge 0} A_1(m,n) x^m q^n \\
Q_2(x) &= \sum_{m,n \ge 0} A_2(m,n) x^m q^n \\
Q_3(x) &= \sum_{m,n \ge 0} A_3(m,n) x^m q^n.
\end{align*}
Upon specializing $x \mapsto 1$ (``forgetting'' the number of parts of the colored partitions), we see that we can rewrite the conclusions of Theorems~\ref{thm:main1} and~\ref{thm:main2} as
\begin{align} Q_1(1) &= \frac{1} {\left(q;q^2\right)_\infty\left(q^1,q^4;q^5\right)_\infty} \label{eq:Q11} \\
Q_3(1) &= \frac{1}{\left(q;q^2\right)_\infty\left(q^2,q^3;q^5\right)_\infty}. \label{eq:Q31}
\end{align}
The three generating functions defined above satisfy the following set of functional equations:
\begin{align} \label{eq:Q1x}
Q_1(x) &=xq^2\,Q_3\left(xq^2\right)
+xq\,Q_3\left(xq\right)
+xq\,Q_2\left(xq\right)
+Q_1\left(xq\right) \\ \label{eq:Q2x}
Q_2(x) &=xq^2\,Q_3\left(xq^2\right)
+xq\,Q_2\left(xq\right)
+Q_1\left(xq\right) \\ \label{eq:Q3x}
Q_3(x) &= xq\,Q_3\left(xq\right)
+Q_1\left(xq\right)
\end{align}
For example, consider a colored partition in $\Gamma$ counted by $Q_1(x)$ (that is, it has no occurrences of $1_\RR$). Exactly one of the following must be true:
\begin{itemize}
\item It contains $2_\RR$. This means the rest of the colored partition cannot contain $1_\GG$, $1_\BB$, $2_\GG$, $2_\BB$, $3_\RR$, $3_\BB$, or $4_\RR$, and so the rest of the colored partition is counted by $Q_1(xq^2)$. Including the initial part $2_\RR$ gives the term $xq^2\, Q_3(xq^2)$ that appears in~\eqref{eq:Q1x}.
\item It contains $1_\BB$. This means the rest of the colored partition cannot contain $1_\GG$, $2_\RR$, or $2_\BB$, and so the rest of the colored partition is counted by $Q_2(xq)$. Including the initial part $1_\BB$ gives the term $xq\, Q_2(xq)$ that appears in~\eqref{eq:Q1x}.
\item It contains $1_\GG$. This means the rest of the colored partition cannot contain $1_\BB$, $2_\RR$, $2_\BB$, or $3_\RR$, and so the rest of the colored partition is is counted by $Q_3(xq)$. Including the initial part $1_\GG$ gives the term $xq\, Q_3(xq)$ that appears in~\eqref{eq:Q1x}.
\item It has no $1_\GG$, $1_\BB$, or $2_\RR$. This is counted by $Q_1(xq)$.
\end{itemize}
Thus, we have justified~\eqref{eq:Q1x}. The justifications of~\eqref{eq:Q2x} and~\eqref{eq:Q3x} are similar.

The Murray-Miller algorithm can be used to ``untangle'' sets of functional equations~\cite{MurrMill}. See  works by Andrews~\cite{And74} for an early application and Takigiku and Tsuchioka~\cite{TTNan} for a recent application of this algorithm to partition identities. In our case, applying the Murray-Miller algorithm to~\eqref{eq:Q1x},~\eqref{eq:Q2x}, and~\eqref{eq:Q3x} gives
\begin{align} \label{eq:MMQ1x}
&Q_1(x) - 
\left(1 + xq + xq^2\right)Q_1(xq) - xq(1-q-xq^3)Q_1(xq^2) -xq^2\left(1-xq^2\right)Q_1(xq^3) \!\!\!\!\!&&=0 \\ \label{eq:MMQ3x}
&Q_3(x)- \left(1+xq+xq^2\right)Q_3(xq)+x^2q^4\,Q_3\left(xq^2\right)- xq^3\left(1-xq^2\right)Q_3\left(xq^3\right) \!\!\!\!\!&&=0
\end{align}
(we will not need an untangled equation for $Q_2(x)$).
\section{Proofs}

Here is an outline of how our proofs will proceed:
\begin{itemize}
\item Show that expressions involving $S_{a,b,c}(x)$ equal certain bivariate generating functions for pairs of (ordinary and uncolored) partitions, by deducing functional equations for these pairs of partitions and showing that the $S_{a,b,c}(x)$  expressions satisfy those functional equations.
\item Show that expressions involving $T_{a,b,c}(x)$ equal $Q_1(x)$ and $Q_3(x)$, by showing that they satisfy~\eqref{eq:MMQ1x} and~\eqref{eq:MMQ3x}, respectively.
\item Set $x=1$ in the expressions above, and conclude that Theorems~\ref{thm:main1} and~\ref{thm:main2} are true.
\end{itemize}
Our proofs showing the expressions satisfy the functional equations all use Maple to write the desired equations as a linear combination of the atomic relations from section 2.1. Verifying this would be impractical to do by hand, but is easy with a computer.

\subsection{Proof of Theorem~\ref{thm:main1}}
We will now prove Theorem~\ref{thm:main1}.

\begin{proof}

First, we will prove:
\begin{align} \label{eq:mainx3}
(1+xq)S_{3,0,1}(x) = \left(-xq;q\right)_\infty\sum_{n\ge0} \frac{x^nq^{n^2}}{(q;q)_n}
\end{align}
We will do this by showing that both sides satisfy the following functional equation:
\begin{align} \label{eq:funceq2x}
f(x) = (1+xq)f(xq)+xq\left(1+xq\right)\left(1+xq^2\right)f\left(xq^2\right)
\end{align}
The right side of~\eqref{eq:mainx3} is the generating function for pairs of partitions $(\lambda,\mu)$, where $\lambda$ is a partition into distinct parts and $\mu$ is a partition where consecutive parts differ by at least 2 (in other words, $\mu$ is counted by the sum side of the first Rogers-Ramanujan identity). To see that the right side of~\eqref{eq:mainx3} satisfies~\eqref{eq:funceq2x}, simply ``peel off'' the smallest parts of these partitions:
\begin{itemize}
\item If the smallest part of $\lambda$ is at least 2 and the smallest part of $\mu$ is at least 2, then the partition pair is counted by $f(xq)$.
\item If the smallest part of $\lambda$ is 1 and the smallest part of $\mu$ is at least 2, then, after deleting the 1 from $\lambda$, the partition pair will be counted by $f(xq)$. Thus, after restoring the 1 to $\lambda$, the original partition pair is counted by $xq\,f(xq)$.
\item If the smallest part of $\lambda$ is at least 2 and the smallest part of $\mu$ is 1, then first delete the 1 from $\mu$. The resulting partition pair will be counted by $f\left(xq\right)$, so after restoring the 1 to $\mu$, the original partition pair is counted by $xq\,f\left(xq^2\right)$.
\item If the smallest part of $\lambda$ is 2 and the smallest part of $\mu$ is 1, then first delete the 2 from $\lambda$ and the 1 from $\mu$. The resulting partition pair will be counted by $f\left(xq^2\right)$, so after restoring the 2 to $\lambda$ and the 1 to $\mu$, the original partition pair is counted by $x^2q^3\,f\left(xq^2\right)$.
\item If the smallest part of $\lambda$ is 1, $\lambda$ does not contain a 2, and the smallest part of $\mu$ is 1, then first delete the 1 from $\lambda$ and the 1 from $\mu$. The resulting partition pair will be counted by $f\left(xq^2\right)$, so after restoring the 1 to $\lambda$ and the 1 to $\mu$, the original partition pair is counted by $x^2q^2\,f\left(xq^2\right)$.
\item If the smallest part of $\lambda$ is 1, $\lambda$ contains a 2, and the smallest part of $\mu$ is 1, then first delete the 1 and the 2 from $\lambda$ and the 1 from $\mu$. The resulting partition pair will be counted by $f\left(xq^2\right)$, so after restoring the 1 and the 2 to $\lambda$ and the 1 to $\mu$, the original partition pair is counted by $x^3q^4\,f\left(xq^2\right)$.
\end{itemize}
Hence, we conclude
\begin{align*}
f(x) &= f(xq) + xq\, f(xq) + xq\,f\left(xq^2\right) + x^2q^3\,f\left(xq^2\right) + x^2q^2\,f\left(xq^2\right) +x^3q^4f\left(xq^2\right) \\
&= (1+xq)f(xq) + xq\left(1 + xq\right)\left(1 + xq^2\right)f\left(xq^2\right),
\end{align*}
as desired.
Substituting $(1+xq)S_{3,0,1}(x)$ into~\eqref{eq:funceq2x} and using~\eqref{eq:Sabcshift} produces
\begin{equation*}
(1+xq)S_{3,0,1} - (1+xq)(1+xq^2)S_{5,1,2}-xq\left(1+xq\right)\left(1+xq^2\right)
(1+xq^3)S_{7,2,3} =0
\end{equation*}
or
\begin{equation*}
S_{3,0,1} - (1+xq^2)S_{5,1,2}-xq\left(1+xq^2\right)
(1+xq^3)S_{7,2,3} =0
\end{equation*}
This can be written as a linear combination of atomic relations. Maple code to verify this is in {\tt MC6\_proof\_A.txt}. Thus, after checking that initial conditions are satisfied (both sides of~\eqref{eq:mainx3} equal 1 when either $x=0$ or $q=0$), we conclude that~\eqref{eq:mainx3} is true.

Next, we will show that $T_{1,0,1}(x)+xq\,T_{3,1,2}(x)$ is the bivariate generating function for the colored partitions discussed in Theorem~\ref{thm:main1}; that is,
\begin{align} \label{eq:37}
T_{1,0,1}(x)+xq\,T_{3,1,2}(x) = Q_1(x).
\end{align}
Substituting $T_{1,0,1}(x)+xq\,T_{3,1,2}(x)$ into~\eqref{eq:MMQ1x} and using~\eqref{eq:Sabcshift} produces
\begin{align*}  \notag
 &T_{1,0,1}
 +xq\,T_{3,1,2}
-\left(1 + xq + xq^2\right)\left(T_{2,1,2}+xq^2\,T_{4,2,3}\right)
-xq\left(1-q-xq^3\right)\left(T_{3,2,3}+xq^3\,T_{5,3,4}\right) \\
&\quad \quad \quad -xq^2\left(1-xq^2\right)\left(T_{4,3,4}+xq^4\,T_{6,4,5}\right)=0
\end{align*}
This can be written as a linear combination of atomic relations. Maple code to verify this is in {\tt MC6\_proof\_B.txt}. Again, after verifying initial conditions are satisfied (both sides of~\eqref{eq:37} equal 1 when either $x=0$ or $q=0$), we conclude that~\eqref{eq:37} is true. Finally, we can prove 
\begin{align} \label{eq:39}
T_{1,0,1}(x)+q\,T_{3,1,2}(x)-(1+q)T_{3,0,1}(x) =0   
\end{align}
by writing it as a linear combination of atomic relations. This can also be verified in Maple (see {\tt MC6\_proof\_E.txt}).

To finish the proof, we use~\eqref{eq:mainx3},~\eqref{eq:37}, and~\eqref{eq:39}, plus~\eqref{eq:29}, Euler's odd-distinct identity, and the first Rogers-Ramanujan identity~\cite{RogRR} to conclude
\begin{align*}
 Q_1(1) &= T_{1,0,1}(1)+q\,T_{3,1,2}(1) \\
 &= (1+q)T_{3,0,1}(1)  \\
 &= (1+q)S_{3,0,1}(1)  \\
&= \left(-q;q\right)_\infty\sum_{n\ge0} \frac{q^{n^2}}{(q;q)_n} \\
&=\frac{1}{\left(q;q^2\right)_\infty\left(q^1,q^4;q^5\right)_\infty},
\end{align*}
thus completing the proof of Theorem~\ref{thm:main1} (see also~\eqref{eq:Q11}).

\end{proof}

\subsection{Proof of Theorem~\ref{thm:main2}}
We will now prove Theorem~\ref{thm:main2}.

\begin{proof}
First, we will prove:
\begin{align} \label{eq:mainx1}
S_{2,0,1}(x) = \left(-xq;q\right)_\infty\sum_{n\ge0} \frac{x^nq^{n^2+n}}{(q;q)_n}
\end{align}
We will do this by showing that both sides satisfy the following functional equation:
\begin{align} \label{eq:funceq1x}
g(x) = (1+xq)g(xq)+xq^2\left(1+xq\right)\left(1+xq^2\right)g\left(xq^2\right)
\end{align}
The right side of~\eqref{eq:mainx1} is the generating function for pairs of partitions $(\lambda,\mu)$, where $\lambda$ is a partition into distinct parts and $\mu$ is a partition where consecutive parts differ by at least 2 and all parts are at least 2 ($\mu$ is counted by the sum side of the second Rogers-Ramanujan identity). To see that the right side of~\eqref{eq:mainx1} satisfies~\eqref{eq:funceq1x}, simply ``peel off'' the smallest parts of these partitions:
\begin{itemize}
\item If the smallest part of $\lambda$ is at least 2 and the smallest part of $\mu$ is at least 3, then the partition pair is counted by $g(xq)$.
\item If the smallest part of $\lambda$ is 1 and the smallest part of $\mu$ is at least 3, then, after deleting the 1 from $\lambda$, the partition pair will be counted by $g(xq)$. Thus, after restoring the 1 to $\lambda$, the original partition pair is counted by $xq\,g(xq)$.
\item If the smallest part of $\lambda$ is at least 3 and the smallest part of $\mu$ is 2, then first delete the 2 from $\mu$. The resulting partition pair will be counted by $g\left(xq^2\right)$, so after restoring the 2 to $\mu$, the original partition pair is counted by $xq^2\,g\left(xq^2\right)$.
\item If the smallest part of $\lambda$ is 2 and the smallest part of $\mu$ is 2, then first delete the 2 from $\lambda$ and the 2 from $\mu$. The resulting partition pair will be counted by $g\left(xq^2\right)$, so after restoring the two 2s to $\lambda$ and $\mu$, the original partition pair is counted by $x^2q^4\,g\left(xq^2\right)$.
\item If the smallest part of $\lambda$ is 1, $\lambda$ does not contain a 2, and the smallest part of $\mu$ is 2, then first delete the 1 from $\lambda$ and the 2 from $\mu$. The resulting partition pair will be counted by $g\left(xq^2\right)$, so after restoring the 1 to $\lambda$ and the 2 to $\mu$, the original partition pair is counted by $x^2q^3\,g\left(xq^2\right)$.
\item If the smallest part of $\lambda$ is 1, $\lambda$ contains a 2, and the smallest part of $\mu$ is 2, then first delete the 1 and the 2 from $\lambda$ and the 2 from $\mu$. The resulting partition pair will be counted by $g\left(xq^2\right)$, so after restoring the 1 and the 2 to $\lambda$ and the 2 to $\mu$, the original partition pair is counted by $x^3q^5\,g\left(xq^2\right)$.
\end{itemize}
Hence, we conclude
\begin{align*}
g(x) &= g(xq) + xq\,g(xq) + xq^2\,g\left(xq^2\right) + x^2q^4\,g\left(xq^2\right) + x^2q^3\,g\left(xq^2\right) +x^3q^5\,g\left(xq^2\right) \\
&= (1+xq)g(xq) + xq^2\left(1 + xq\right)\left(1 + xq^2\right)g\left(xq^2\right),
\end{align*}
as desired. Substituting $S_{2,0,1}(x)$ into~\eqref{eq:funceq1x} and using~\eqref{eq:Sabcshift} produces
\begin{equation*}
S_{2,0,1}-(1+xq)S_{4,1,2}-xq^2(1+xq)\left(1+xq^2\right)S_{6,2,3}=0.
\end{equation*}
This can be written as a linear combination of atomic relations. Maple code to verify this is in {\tt MC6\_proof\_C.txt}. Since the initial conditions are satisfied (both sides of~\eqref{eq:mainx1} equal 1 when either $x=0$ or $q=0$), we conclude that~\eqref{eq:mainx1} is true.

Now, we will show that $T_{2,0,1}(x)$ is the bivariate generating function for the colored partitions discussed in Theorem~\ref{thm:main2}; that is, 
\begin{equation}   \label{eq:mainx2} T_{2,0,1}(x) = Q_3(x).\end{equation}
Substituting $T_{2,0,1}(x)$ into~\eqref{eq:MMQ3x} and using~\eqref{eq:Tabcshift} produces
\begin{equation*}
T_{2,0,1} - (1+xq+xq^2)T_{3,1,2}+x^2q^4\,T_{4,2,3} -xq^3(1-xq^2)T_{5,3,4} = 0.
\end{equation*}
This can be written as a linear combination of atomic relations. Unlike the others in this paper, the linear combination is small enough to type out (and the truly motivated reader could verify by hand):
\begin{align*}
&\frac{xq}{q^2 + 1}\relh^1_{2, 0, 1}
-xq\,\relh^1_{2, 0, 3}
+\frac{xq^3}{(q - 1)(q^2 + 1)} \relh^2_{4, 0, 1}
-\frac{xq^4}{(q - 1)(q^2 + 1)} \relh^2_{4, 0, 3}
-\frac{xq(q^2 - q + 1)}{(q - 1)(q^2 + 1)} \relh^3_{4, 0, 1} \\
&+\frac{xq^3}{(q - 1)(q^2 + 1)} \relh^3_{4, 2, 1} -\frac{(q^5 - q^4 + 2q^3 - 2q^2 - qx + q - 1)}{(q - 1)(q^2 + 1)^2}\relh^6_{2, 0, 1} \\
&+\frac{(q^4 - q^3 + 2q^2 - 2q + 1)xq^2}{(q - 1)(q^2 + 1)^2}
\relh^6_{3, 1, 2}-\frac{xq^2}{(q - 1)(q^2 + 1)^2} \relh^8_{2, 0, 1}- \relh^8_{2, 0, 3}
-\frac{xq^2}{(q - 1)(q^2 + 1)^2} \relh^8_{3, 1, 2} \\
&-\frac{xq^3}{(q^2 + 1)^2} \relh^9_{2, 0, 1}
+(1+xq) \relh^9_{2, 0, 3}
-\frac{xq^3}{(q^2 + 1)^2}  \relh^9_{3, 1, 2}\\
&= T_{2,0,1} - (1+xq+xq^2)T_{3,1,2}+x^2q^4T_{4,2,3} -xq^3(1-xq^2)T_{5,3,4}. 
\end{align*}
In any case, Maple code to verify this is in {\tt MC6\_proof\_D.txt}. We can conclude~\eqref{eq:mainx2} is true after checking the initial conditions (both sides of~\eqref{eq:mainx2}equal 1 when either $x=0$ or $q=0$).

To finish off the proof, substituting $x=1$ in~\eqref{eq:mainx1} and applying Euler's odd-distinct identity and the second Rogers-Ramanujan identity~\cite{RogRR} gives
\begin{align*}
\sum_{i,j,k\ge 0}
\frac{q^{\frac 32i^2+\frac12j^2+\frac12k^2+ij+ik+jk-\frac12i+\frac12j+\frac12k+2i+k}}{\left(q;q\right)_{i+j+k}}\qtri{i+j+k}ijk{q^2} = \frac{1}{\left(q;q^2\right)_\infty\left(q^2,q^3;q^5\right)_\infty}.
\end{align*}
Meanwhile, substituting $x=1$ in~\eqref{eq:mainx2} gives
\begin{align*}
\sum_{i,j,k\ge 0}
\frac{q^{\frac 32i^2+\frac12j^2+\frac12k^2+ij+ik+jk-\frac12i+\frac12j+\frac12k+2i+k}}{\left(q;q\right)_{i+j+k}}\qtri{i+j+k}ijk{q^2} = \sum_{n\ge0} A^\ast(n) q^n = Q_3(1).
\end{align*}
We conclude that 
\begin{align*}
Q_3(1)=\sum_{n\ge0} A^\ast(n) q^n = \frac{1}{\left(q;q^2\right)_\infty\left(q^2,q^3;q^5\right)_\infty},
\end{align*}
and we have accordingly completed our proof of Theorem~\ref{thm:main2} (see also~\eqref{eq:Q31}).
\end{proof}

\section{Discussion}

Computations suggest the following conjectured variation of Theorem~\ref{thm:main1}:
\begin{conj}
Let $A^\ast_1(m,n)$ be the number of three-colored partitions in $\Gamma$ of $n$ that have no occurrences of $1_\RR$, where $m$ is the number of parts colored green or blue plus twice the number of parts colored red. Define 
\begin{align*}
Q^\ast_1(y)= \sum_{m,n \ge 0} A^\ast_1(m,n) y^m q^n.
\end{align*}
Then, we conjecture $
Q^\ast_1(y) = (1+yq)S_{3,0,1}(y).$
\end{conj}

Perhaps equally interesting as the identities themselves is the identity-proving paradigm suggested by this paper and the previous paper of K. Baker, Kanade, C. Sadowski, and the author~\cite{BKRS}.  Suppose that one wants to prove the identity 
\begin{align*}
\sum_{i,j,k\ge 0}
\frac{q^{\frac 32i^2+\frac12j^2+\frac12k^2+ij+ik+jk-\frac12i+\frac12j+\frac12k+\left(2i+k\right)}}{\left(q;q\right)_{i+j+k}}\qtri{i+j+k}ijk{q^2}
=
\frac{1}{\left(q;q^2\right)_\infty\left(q^2,q^3;q^5\right)_\infty}
\end{align*}
that was conjectured through computer experimentation. One approach (roughly the one followed by this paper) is to insert in one auxiliary variable for each summation variable on the left (here, $u,$ $v$, and $w$), and rewrite the right side in a way that is amenable to the introduction of auxiliary variables (here, $y$ and $z$). This produces
\begin{align*}
\sum_{i,j,k\ge 0}
\frac{q^{\frac 32i^2+\frac12j^2+\frac12k^2+ij+ik+jk-\frac12i+\frac12j+\frac12k+\left(2i+k\right)}}{\left(q;q\right)_{i+j+k}}\qtri{i+j+k}ijk{q^2} u^iv^jw^k
=
\left(-yq;q\right)_\infty\sum_{n\ge0} \frac{z^nq^{n^2}}{(q;q)_n}
\end{align*}
Then, the goal is to experimentally find a way to specialize $u,$ $v$, $w,$ $y$, and $z$ to make the above identity true. There is no guarantee that this process will work, but if it does (for example, in this case, by setting $u=x^2$ and $v=x$, $w=x$, $y=x$, and $z=x$), then one can deduce a functional equation (from combinatorial principles) to govern the resulting bivariate right side (in our case,~\eqref{eq:funceq1x}) and then use the method of atomic relations to show that the multisum satisfies the same functional equation.

There is currently a large amount of interest in proving the modularity of families of Nahm sums (and slight variants of Nahm sums). See, for instance, works by Z. Cao, S.-P. Cui, J. Q. D. Du, Z. Li, Y. Mizuno, H. Rosengren, C. Shi, B. Wang, and L. Wang~\cite{CRW,LiWa,Mizuno,WW_Miz,CD,SW,CaoWang}, among others. Most of the identities in these papers have been proved, but there are some conjectures left open. Perhaps the method outlined here could be adapted to prove of future conjectures in this vein.

\providecommand{\oldpreprint}[2]{\textsf{arXiv:\mbox{#2}/#1}}\providecommand{\preprint}[2]{\textsf{arXiv:#1 [\mbox{#2}]}}

\end{document}